%% file: main.tex
\documentclass[12pt]{article}
\usepackage{geometry}
\usepackage{graphicx}
\usepackage{amsmath, amsfonts, amsthm, amssymb}
\usepackage{url}
\usepackage{tikz} 
\usepackage[colorlinks=true, allcolors=blue]{hyperref} 
\usepackage{aliascnt}
\usepackage[colorlinks=true, allcolors=blue]{hyperref}
\usepackage[nameinlink]{cleveref}
\theoremstyle{definition}

\newtheorem{theorem}{Theorem}

\newaliascnt{lemma}{theorem}
\newtheorem{lemma}[lemma]{Lemma}
\aliascntresetthe{lemma}

\newaliascnt{corollary}{theorem}
\newtheorem{corollary}[corollary]{Corollary}
\aliascntresetthe{corollary}

\newaliascnt{conjecture}{theorem}
\newtheorem{conjecture}[conjecture]{Conjecture}
\aliascntresetthe{conjecture}

\newaliascnt{proposition}{theorem}
\newtheorem{proposition}[proposition]{Proposition}
\aliascntresetthe{proposition}

\newaliascnt{remark}{theorem}
\newtheorem{remark}[remark]{Remark}
\aliascntresetthe{remark}

\newaliascnt{definition}{theorem}
\newtheorem{definition}[definition]{Definition}
\aliascntresetthe{definition}

\newaliascnt{example}{theorem}
\newtheorem{example}[example]{Example}
\aliascntresetthe{example}

\crefname{theorem}{theorem}{theorems}
\Crefname{theorem}{Theorem}{Theorems}
\crefname{lemma}{lemma}{lemmas}
\Crefname{lemma}{Lemma}{Lemmas}
\crefname{corollary}{corollary}{corollaries}
\Crefname{corollary}{Corollary}{Corollaries}
\crefname{definition}{definition}{definitions}
\Crefname{definition}{Definition}{Definitions}
\crefname{proposition}{proposition}{propositions}
\Crefname{proposition}{Proposition}{Propositions}
\crefname{remark}{remark}{remarks}
\Crefname{remark}{Remark}{Remarks}
\crefname{example}{example}{examples}
\Crefname{example}{Example}{Examples}
\crefname{conjecture}{conjecture}{conjectures}
\Crefname{conjecture}{Conjecture}{Conjectures}

\usepackage{enumitem}
\makeatletter
\def\saveenum{\xdef\@savedenum{\the\c@enumi\relax}}
\def\resetenum{\global\c@enumi\@savedenum}
\makeatother

\usepackage[dvipsnames]{xcolor} 
\usepackage{comment}

\usepackage{framed} 
\usepackage{subcaption} 
\usepackage{booktabs} 
\usepackage{arydshln} 
\usepackage{mathtools} 

\newcommand{\demph}[1]{\emph{{\color{black}#1}}}

\newcommand{\QQ}{\mathbb{Q}}

\newcommand{\stargraph}{\mathcal{S}}

\newcommand{\energy}{\mathcal{E}}
\newcommand{\Spec}{\text{Spec}}

\DeclareMathOperator{\Gal}{Gal}
\DeclareMathOperator{\tr}{tr}

\title{\bf Galois Actions on Energy-Harmonic Functions}

\author{Bella Finkel}

\begin{document}
\date{}
\maketitle

\begin{abstract}
    We study the action of $\Gal(\Bar{\QQ}/\QQ)$ on the solutions to the discrete Dirichlet problem for a graph with fixed edge energies. We establish transitivity of the Galois action outside a thin set in the parameter space and demonstrate that the exceptional thin set may be Zariski dense, so the conclusion cannot in general be replaced by one that holds on a Zariski-open set. 
    Identifying the parameter locus where the Galois action is transitive constrains the possibilities for linking combinatorial and arithmetic properties of the discrete energy-harmonic Laplace equations and other systems of the same form.
\end{abstract}
\section{Introduction}
Determining the Galois or monodromy group of a family of equations reveals structural features of the solutions. It admits study via numerical algebraic geometry \cite{HauensteinRodriguezSottile,leykin2009galois} and itself bounds the complexity of finding solutions \cite{Duff2025}.
Geometrically, the question is which permutations of the solutions arise from varying the parameters along closed paths that avoid the branch locus. Algebraically, Harris \cite{Harris1979} showed that the monodromy group is isomorphic to the Galois group of the Galois closure of the corresponding extension of complex function fields for appropriate parameter and solution spaces. 
If the parameter space is rational, then the Hilbert irreducibility theorem implies that the Galois group of the Galois closure of the smallest field of definition of a solution is 
precisely the generic arithmetic Galois group with positive probability.
When a  
system of equations 
has a combinatorial interpretation, it is natural to ask whether there is a Galois-invariant 
combinatorial property that tracks meaningful features of the problem.

The most basic instance of this question is to ask for a combinatorial invariant that gives the cardinality of its solution set. For a graph $G$ with distinct rational boundary values and positive rational edge weights, the number of solutions to the
Dirichlet problem on $G$
equals the number of \emph{compatible orientations} of $G$ in the sense of \cite[Section 1.2]{KenyonEnharmonic}--indeed, each solution induces a compatible orientation \cite[Theorem 1]{KenyonEnharmonic}. In this sense, the compatible orientations  
of $G$ furnish a set of combinatorial objects permuted by the absolute Galois group $\Gal(\Bar{\QQ}/\QQ)$.

One piece of evidence that the discrete Dirichlet problem is arithmetically rich is that the induced action of the Galois group of the totally real algebraic numbers over $\mathbb{Q}$ is faithful as the graph and parameters vary: every totally real number field is generated by the solutions to the Dirichlet problem for some graph \cite[Theorem 6]{KenyonEnharmonic}.
Conjecture 1 of \cite{KenyonEnharmonic} proposes that, for a 3-connected graph with boundary, this action is transitive for every parameter choice 
in a suitable Zariski-open subset of the parameter space.
Our main contribution is a cautionary result in this direction. We provide counterexamples to \cite[Conjecture 1]{KenyonEnharmonic} (\Cref{sec:Zariski–DenseExceptionalLoci}) but demonstrate that Hilbert irreducibility still guarantees transitivity outside of a thin set in the parameter space (\Cref{thm:Conjecture2}). We provide examples characterizing the generic Galois groups and the thin sets (\Cref{sec:GenericGaloisGroups}).

\subsection{Transitive Galois action}

In the case of the energy-harmonic discrete Dirichlet problem, there is an opportunity to identify combinatorial properties of the graph and geometric properties of the parameters with arithmetic properties of the solution set. Other works \cite{gross2014maximum,Early_Pfister_Sturmfels_2025,betti2025graphical, agostini2023likelihood} have investigated how the number of solutions to systems of the same form changes when certain parameters are taken to zero. Even when the cardinality of the solution set is fixed, the choice of parameters can alter its arithmetic behavior.
The purpose of this brief article is to prove the following result establishing transitivity of the Galois action away from a thin set in the parameter space and qualify its potential generalizations.
\begin{theorem}\label{thm:Conjecture2}
    Let $Y$ be the irreducible incidence variety for the discrete enharmonic Laplace equations of $G$.
    There exists a thin set $A \subset \QQ^{|B|+|E|}$ such that 
    for every choice of parameters $$t:=(b_u,\energy_{vw})\in(\QQ^{|B|+|E|})\setminus A,$$ the fiber $Y_t$ of the projection $\pi:Y\to\mathcal{B}$ over $t$ is integral over $\QQ$.
    Hence the absolute Galois group $\Gal(\Bar{\QQ}/\QQ)$ acts transitively on the solutions of the discrete enharmonic Laplace equations with parameters $t$.

\end{theorem}

\subsection{Outline of the paper}
This paper is structured as follows. In \Cref{Sec:EnharmonicIncidenceVarieties}, we describe the discrete Dirichlet problem 
we consider, defined by the discrete enharmonic Laplace equations with rational parameters associated to a graph, and introduce the relevant incidence variety. 
\Cref{sec:HilbertIrreducibilityAndTransitiveGaloisActions} reviews thin sets and establishes how they can be used to appropriately constrain the locus of the incidence variety where the Galois action is transitive, culminating in the proof of \Cref{thm:Conjecture2}. In \Cref{sec:Zariski–DenseExceptionalLoci}, we show that ``thin" cannot be replaced with ``Zariski open," providing a disproof of \cite[Conjecture 1]{KenyonEnharmonic} in contrast with the positive statement of \Cref{cor:GaloisAction}.
In \Cref{sec:GenericGaloisGroups}, we investigate the Galois groups and illustrate \Cref{thm:Conjecture2} with an example for the star graph with four leaves. \Cref{sec:Conclusion} concludes with a discussion of future directions focused on extensions of our results in applications.

\section{Enharmonic incidence varieties}\label{Sec:EnharmonicIncidenceVarieties}

\subsection{Graphs with boundary and enharmonic functions}\label{sec:Graphs with boundary and enharmonic functions}

In this paper, a \demph{graph} $G=(V,E,B)$ indicates a weighted graph with boundary that is finite, connected, and undirected with no loops or multiple edges. 
We say that $G$ has vertex set $V$, edge set $E$ and boundary $B$, and assume that the designated set of boundary vertices $B\subset V$ satisfies $|B|\geq2$ and that each edge of $G$ lies on a simple path from some vertex in $B$ to another vertex in $B$. 
Each edge $vw:=\{v,w\}$ is assigned a nonzero rational weight $\energy_{vw}$ called the \demph{(edge) energy}.
A \demph{boundary value function} $b=(b_u)_{u\in B}:B\to\QQ$ is an injective function that assigns a rational number to each boundary vertex. We call the energies and boundary values the \emph{parameters} of $G$.

Corresponding to such a graph is a discrete Dirichlet problem where one seeks a function 
$h:V\to\mathbb{R}$ 
harmonic on $V\setminus B$ such that $h|_{B}=b.$ 
If one fixes the desired edge energies of the harmonic function, the problem becomes to determine
an \demph{enharmonic} (``energy-harmonic'') function $h$. For a finite connected graph, this is defined as a function 
with $h|_B=b$ that satisfies the \demph{discrete enharmonic Laplace equations}
\begin{equation}\label{eq:Lh}
    \sum_{w\sim v}\frac{\mathcal{E}_{vw}}{h(v)-h(w)}=0
\end{equation}
at every interior vertex $v$, where the sum is taken over all vertices $w$ adjacent to $v$. The equations \eqref{eq:Lh} are algebraic, so when the energies are fixed the function $h$ is permitted to be complex-valued.
For an expression in terms of \emph{edge conductances} equivalent to \eqref{eq:Lh} more familiar from the perspective of finding 
the unique function minimizing the Dirichlet energy, see \cite[Section 1.2]{KenyonEnharmonic}. If every vertex is a boundary vertex, then $G$ does not have an associated Dirichlet problem.

\subsection{Parameter space and incidence variety}

An incidence variety for the Dirichlet problem on a graph with fixed edge energies can be constructed as follows. For each vertex $v\in V$, define 
\begin{equation*}
    H_v =
\begin{cases}
x_v & \text{if } v \in V \setminus B, \\
b_v & \text{if } v \in B
\end{cases}
\end{equation*}
and
\begin{equation*}
    E^\text{int}=\big\{vw\in E:\{v,w\}\not\subseteq B\big\}.
\end{equation*}
Give $\mathbb{A}_{\QQ}^{|V|+|E|}$ coordinates
\[
(x_v)_{v\in V\setminus B},\qquad
(b_u)_{u\in B},\qquad \text{and} \qquad
(\mathcal E_{vw})_{vw\in E}.
\]
Define the enharmonic incidence variety $Y$  to be the 
locally closed subvariety
\begin{equation}\label{eq:IncidenceVariety}
\begin{aligned}
    Y:=\mathcal{V}\bigg(\sum_{w\sim v}\energy_{vw}\prod_{\substack{w' \sim v\\ w' \neq w}}(H_v-H_{w'}): v\in V\setminus B\bigg)
    \bigcap
    \\
    D\bigg(\prod_{\{v,w\}\in E^{\text{int}}}(H_v-H_w)\prod_{\{u,u'\}\in \binom{B}{2}}(b_u-b_{u'})\bigg)\subset\mathbb{A}_{\QQ}^{|V|+|E|}
\end{aligned}
\end{equation}
where $D(f)$ denotes the distinguished open subset on which $f$ does not vanish. We write $\QQ[Y]$ to mean the ring of regular functions on $Y$ as an affine variety.
In the product over $E^{\text{int}}$, each unoriented edge appears only once.
\begin{lemma}\label{lem:IrreducibleIncidenceVariety}
    The incidence variety $Y$ for the discrete enharmonic Laplace equations with pairwise distinct boundary values is irreducible for any $G$ as defined in \Cref{sec:Graphs with boundary and enharmonic functions}. 
\end{lemma}
\begin{proof}
    After a change of variables 
    \begin{equation*}
        \omega_{vw}=\frac{\energy_{vw}}{h(v)-h(w)}
    \end{equation*}
    to edge flows $\omega$ as in \cite[proof of Theorem 1]{KenyonEnharmonic}, 
    \eqref{eq:Lh} becomes a system of homogeneous linear equations 
    with constant coefficients. Therefore $Y$ is a trivial vector bundle over the open subset of vertex values defined by the nonvanishing conditions in \eqref{eq:IncidenceVariety}. Since this subset is irreducible, $Y$ is also.
\end{proof}

For a statement like \Cref{thm:Conjecture2} to make sense, the generic fiber in $Y$ over the parameter space must be finite. This follows from a well-known theorem of Orlik and Terao \cite[Theorem 1.1]{OrlikTerao} once it has been established that the enharmonic functions are the critical points of the master function of a Dirichlet arrangement \cite[Theorem 1.4]{lutz2019electrical}. When the edge energies are strictly positive, Varchenko's theorem \cite[Theorem 1.2.1]{Varchenko1995} 
is sufficient, as was established independently by \cite[Theorem 1]{KenyonEnharmonic} and \cite[Corollary 4]{KenyonEnharmonic}.

\section{Hilbert irreducibility and transitive Galois actions}\label{sec:HilbertIrreducibilityAndTransitiveGaloisActions}

\subsection{Thin sets}

Thin sets were defined by Serre \cite[Definition 3.1.1]{serre} in his geometric formulation of Hilbert's irreducibility theorem to describe exceptional loci where the Galois group of a cover may fail to be preserved as the Galois group of a fiber.
In the present article, a salient aspect of thinness is that it provides a notion of membership in a set being ``unlikely" that encompasses certain Zariski-dense subsets.
In conjunction with the examples in \Cref{sec:Zariski–DenseExceptionalLoci}, the phrasing of \Cref{thm:Conjecture2} in terms of thin sets provides a useful contrast to the genericity conclusion of \cite[Conjecture 1]{KenyonEnharmonic}. In \Cref{sec:GenericGaloisGroups}, we will characterize a thin set for \Cref{thm:Conjecture2} when $G$ is a star graph with four leaves.

\begin{definition}{\cite[Definition 3.1.1]{serre}}\label{def:ThinSet}   Let $K$ be a field of characteristic zero, $\mathcal{Z}$ an irreducible algebraic variety over $K$, and $\mathcal{Z}(K)$ the set of $K$-rational points of $\mathcal{Z}$. Suppose $A$ is a subset of $\mathcal{Z}(K)$. If $A$ is not Zariski-dense in $\mathcal{Z}$, then $A$ is said to be a \demph{thin set of type I}. Otherwise, if there is an irreducible variety $\mathcal{Z}'$ of the same dimension as $\mathcal{Z}$ and a generically surjective morphism $\pi:\mathcal{Z}'\to \mathcal{Z}$ of degree at least two with $A\subset \pi(\mathcal{Z}'(K))$, then $A$ is said to be a \demph{thin set of type II}. 

The subset $A$ is called simply a \demph{thin set} if it is contained in a finite union of thin sets of types I and II.
\end{definition}

\subsection{The generic incidence extension}

By a \demph{branched cover}, we mean a dominant separable map between irreducible varieties of the same dimension; see \cite[Section 11.1]{EisenbudHarris2025}.

\begin{lemma}\label{lem:branchedcover}
    Set 
    \[\mathcal{B}:=D\big(\prod_{\{u,u'\}\in \binom{B}{2}}(b_u-b_{u'})\big)\subset\mathbb{A}_{\QQ}^{|B|+|E|}\]
    and let $Y$ be the incidence variety defined in \eqref{eq:IncidenceVariety}.
    Then the projection map $\pi:Y\to \mathcal{B}$ is a branched cover of finite degree. Furthermore, there is an affine open subset $W^\circ\subseteq \mathcal{B}$ such that
    \begin{equation}\label{eq:PrimitiveElementEquation}
        \QQ[Y^\circ]=\QQ[W^\circ][X]/\langle f\rangle
    \end{equation}
    where $Y^\circ=\pi^{-1}(W^\circ)$ and  $f\in\QQ[W^\circ][X]$ is a monic polynomial irreducible over $\QQ(W^\circ)$. If $\pi$ has degree $d$, then so does $f$.
    \end{lemma}

\begin{proof}
      The geometric points of the fiber of $\pi$ over a point in $\mathcal{B}(\mathbb{Q})$ consist exactly of the enharmonic functions solving \eqref{eq:Lh} for that tuple of parameters. The number $d$ of enharmonic functions of a fixed graph with generic parameters is 
      nonzero by \cite[Proposition 3.15]{lutz2019electrical}, so the generic fiber of  $\pi$ is nonempty. Thus $\pi(Y)$ is Zariski dense in $\mathcal{B}$, so $\pi$ is dominant. 
      By \cite[Corollary 1.5]{lutz2019electrical} it is generically finite, so $Y$ has dimension $|B|+|E|$. Since $\mathcal{B}$ and $Y$ are irreducible, $\pi$ is a degree-$d$ branched cover.

     The function field $\QQ(Y)$ is a finite extension of the function field $\QQ(\mathcal{B})$ of degree $d$. 
    The primitive element theorem guarantees there is an element $\alpha\in\QQ(Y)$ such that $\QQ(Y)=\QQ(\mathcal{B})[X]/(\mathfrak{m}_\alpha)$ where 
    \[
    \mathfrak{m}_\alpha(X)=X^d+a_1X^{d-1}+\cdots+a_d,\quad a_i\in\mathbb{Q}(\mathcal{B})
    \]
    is the minimal polynomial of $\alpha$. 

    One can find an affine open subset $W^\circ\subseteq\mathcal{B}$ such that $\pi|_{Y^\circ}$ is finite and $\alpha$ is regular on $\pi^{-1}(W^\circ)$ and integral over $\QQ[W^\circ].$ If $\QQ[Y^\circ]$ and $\QQ[W^\circ]$ are the coordinate rings, there is a monic irreducible polynomial $f\in\QQ[W^\circ][X]$ of degree $d$ such that 
    \begin{equation}
        \QQ[Y^\circ]=\QQ[W^\circ][X]/\langle f\rangle.
    \end{equation}
\end{proof}

\subsection{Proof of a transitive Galois action outside of a thin set}

To understand whether the Galois action is transitive for particular energies and boundary values, we need to know under what conditions the fiber over a rational point in parameter space is integral over $\QQ$. To this end, we recall the following definition and proposition from \cite[Chapter 3]{serre}.
\begin{definition}
    Let $K$ be a field and suppose $\mathcal{Z}$ is a $K$-irreducible variety with function field $K(\mathcal{Z})$. If
\[f(X)=X^n+a_1X^{n-1}+\cdots+a_n,\quad a_i\in K(\mathcal{Z})\]
is an irreducible polynomial over $K(\mathcal{Z})$ and $t\in \mathcal{Z}(K)$ such that $t$ is not a pole of any $a_i$, then $a_i(t)\in K$. We call 
\[f_t(X)=X^n+a_1(t)X^{n-1}+\cdots+a_n(t)\]
the \demph{specialization of $f$ at $t$}.
\end{definition}

\begin{proposition}{\cite[Proposition 3.3.5]{serre}}\label{prop:3.3.5}
    Let $\mathcal{Z}$ be a $\QQ$-irreducible variety with function field $\QQ(\mathcal{Z})$,
    \begin{equation}\label{eq:f(X)}
        f(X)=X^n+a_1X^{n-1}+\cdots+a_n,\quad a_i\in\QQ(\mathcal{Z})
    \end{equation}
     an irreducible polynomial over $\QQ(\mathcal{Z})$, and $H\subset S_n$ the Galois group of $f$. Then there exists a thin set $A\subset \mathcal{Z}(\QQ)$ such that $t\notin A$ implies all of the following:
     \begin{enumerate}[label=\textit{\roman*}.]
         \item no $a_i$ has a pole at $t$,
         \item $f_t(X)$ is irreducible over $\QQ$, and
         \item the specialization $f_t$ has Galois group $H$.
     \end{enumerate}
\end{proposition}

\begin{proof}[Proof of \Cref{thm:Conjecture2}]
Let $W^\circ$ be as in \Cref{lem:branchedcover}. If $t\in W^\circ(\QQ)$ and
    \begin{equation}
        f(X)=X^d+a_1X^{d-1}+\cdots+a_d,\quad a_i\in\mathbb{Q}[W^\circ],
    \end{equation}
    then $a_i(t)\in\mathbb{Q}$ and the specialization of $f$ at $t$ is 
    $f_t(X)=X^d+a_1(t)X^{d-1}+\cdots+a_d(t).$
    The fiber $Y^\circ_t$ of $\pi$ over $t$ is
    \begin{equation}
        Y_t^\circ=\Spec(\kappa(t))\times_{W^\circ}Y^\circ,
    \end{equation}
    where $\kappa(t)$ is the residue field of $W^\circ$ at $t$. Since $t$ is a $\QQ$-rational point, $\kappa(t)=\QQ$. 
    Writing 
    \begin{equation}
        Y^\circ=\Spec(\QQ[W^\circ][X]/\langle f\rangle),
    \end{equation}
     the fiber can be expressed as
    \begin{equation}
        \begin{aligned}
            Y^\circ_t=\Spec\left(\QQ\otimes_{\QQ[W^\circ]}(\QQ[W^\circ][X]/\langle f\rangle)\right).
        \end{aligned}
    \end{equation}
    
    Its coordinate ring is obtained by the base change from $\QQ[W^\circ][X]/\langle f\rangle$ along the evaluation-at-$t$ map $\QQ[W^\circ]\to\QQ$. Hence, for $t\in W^\circ(\QQ)$,
    \begin{equation}\label{eq:CoordinateRingOfTheFiber}
        \QQ[Y^\circ_t]=\QQ[X]/\langle f_t(X)\rangle. 
    \end{equation}
    
    Under the assumption that the parameter space is a $\QQ$-irreducible variety, there exists a thin set $$A'\subset W^\circ(\QQ)$$ such that if $$t\in W^\circ(\QQ)\setminus A',$$ then $f_t(X)$ is irreducible over $\QQ$ \cite[Proposition 3.3.5]{serre}. The complement of $W^\circ(\QQ)$ in $\QQ^{|B|+|E|}$ is thin of type I because it is contained in a proper Zariski closed subset. Therefore
    \begin{equation*}
        A=A'\cup\big(\QQ^{|B|+|E|}\setminus W^\circ(\QQ)\big)
    \end{equation*}
    is a thin set contained in $\QQ^{|B|+|E|}$. 

    Suppose $t\in W^\circ(\QQ)\setminus A$. Then $f_t(X)$ is irreducible over $\QQ$, and \eqref{eq:CoordinateRingOfTheFiber} says that the geometric points of the fiber--that is, the solutions to the discrete enharmonic Laplace equations with parameters $t$--correspond to the roots of $f_t(X)$. Hence $\Gal(\Bar{\QQ}/\QQ)$ acts transitively on the enharmonic functions.
\end{proof}

\subsection{The totally real specialization}
Suppose that the energies and boundary values of $G$ are rational and the energies are strictly positive. 
Denoting the totally real algebraic numbers by $\QQ^{\tr}$, it follows from Varchenko's theorem \cite[Theorem 1.2.1]{Varchenko1995} that the Galois group $\Gal(\QQ^{\tr}/\QQ)$ permutes the enharmonic functions of $G$ \cite[Theorem A.1]{lutz2019electrical}. This fact was proven independently in \cite[Theorem 1 and Corollary 5]{KenyonEnharmonic}, where it was furthermore shown that every totally real number field is generated by the enharmonic functions of some such graph \cite[Theorem 6]{KenyonEnharmonic}. 

\begin{corollary}\label{cor:GaloisAction}
    Let the incidence variety $Y$ be as in \Cref{thm:Conjecture2}. If 
    $$t\in\big(\{b\in\QQ^B:b_u\neq b_v\text{ for } u\neq v\}\times\QQ_{>0}^E\big)\setminus A$$ 
    so that the edge energies are positive, then $\Gal(\QQ^{\tr}/\QQ)$ acts transitively on the solutions of the discrete enharmonic Laplace equations.
\end{corollary}
\begin{proof}[Proof of \Cref{cor:GaloisAction}]
Let $A$ be the thin set from \Cref{thm:Conjecture2}. If 
$$t\in\big(\{b\in\QQ^B:b_u\neq b_v\text{ for } u\neq v\}\times\QQ_{>0}^E\big)\setminus A,$$ 
then \Cref{thm:Conjecture2} implies that the absolute Galois group $\Gal(\Bar{\QQ}/\QQ)$ acts transitively on the solutions to the discrete enharmonic Laplace equations with parameters $t$. Since the edge energies are strictly positive, \cite[Corollary 5]{KenyonEnharmonic} says that every solution takes values in $\mathbb{Q}^{\tr}$. It follows that the action of $\Gal(\overline{\QQ}/\QQ)$ on the solution set factors through $\Gal(\QQ^{\tr}/\QQ)$. The induced action remains transitive.
\end{proof}

\section{Zariski-dense exceptional loci}\label{sec:Zariski–DenseExceptionalLoci}

So far, we have not imposed any conditions on $G$ except to ensure that it corresponds to a nontrivial discrete Dirichlet problem with prescribed nonzero edge energies. A natural next question is whether restricting to graphs with certain properties yields exceptional arithmetic behavior. For instance, is there a hypothesis on $G$ that guarantees the thin set of \Cref{thm:Conjecture2} can be taken as a thin set of type I? One plausible candidate is to require $G$ to be \emph{3-connected}. 

For a graph without boundary, the usual notion of 3-connectedness is the following.
We use the convention that a graph is \demph{connected} if it has exactly one connected component.

\begin{definition}\label{def:3-connected}
    A graph $\Gamma=(V,E)$ is said to be \demph{3-connected if $|V|>3$ and $\Gamma-X$} is connected for every set $X\subset V$ such that $|X|<3$ (see, e.g., \cite[Section 1.4]{diestel2005graph}). 
\end{definition}
For the graphs with boundary considered in this paper, we use the following definition of 3-connectedness, following \cite{KenyonEnharmonic}.
\begin{definition}\label{def:3-connectedWBoundary}
    We say that a graph with boundary $G = (V,E,B)$ is \demph{$3$-connected} if the augmented graph 
\[ G^{B} := \left(V, E\cup \bigl\{\{v,w\}: v,w\in B,\ v\neq w\bigr\} \right) \] 
obtained by adding an edge between every pair of distinct boundary vertices is $3$-connected in the sense of \Cref{def:3-connected}. 
Equivalently, $G^{B}$ remains connected after the removal of any one or two vertices of $V$, whether boundary or interior.
\end{definition}

In general, however, \Cref{thm:Conjecture2} cannot be generalized to guarantee transitivity for parameters in a Zariski-open set of $\mathcal{B}$ even for such graphs. In the remainder of this section, we contrast \Cref{cor:GaloisAction} with the following conjecture. 

\begin{conjecture}\label{conj:Conjecture1}
    \cite[Conjecture 1]{KenyonEnharmonic} If $G$ is 3-connected, the action of $\Gal(\QQ^{\tr}/\QQ)$ on the enharmonic functions is transitive for generic rational boundary data and generic positive rational energies. That is, there exists a Zariski-open set $U\subseteq \mathcal{B}$ such that the action of $\Gal(\QQ^{\tr}/\QQ)$ on the enharmonic functions is transitive for every 
    $$t\in U(\QQ)\cap (\QQ^B\times\QQ_{>0}^E).$$
\end{conjecture}

The following counterexamples demonstrate an obstruction to finding a combinatorial invariant of the graph that determines the arithmetic properties of the solution set for a Zariski-open subset of the parameter space.
\begin{definition}
    The \demph{star graph} $\stargraph_{n}$ is the $(n-1)$-leafed tree with a single interior vertex whose boundary is precisely the set of leaves.
\end{definition}

\begin{figure}
    \centering
    \input{S8}
    \caption{Star graph on eight vertices. Boundary vertices are marked in blue and interior vertices are marked in purple.}
    \label{fig:S8}
\end{figure}
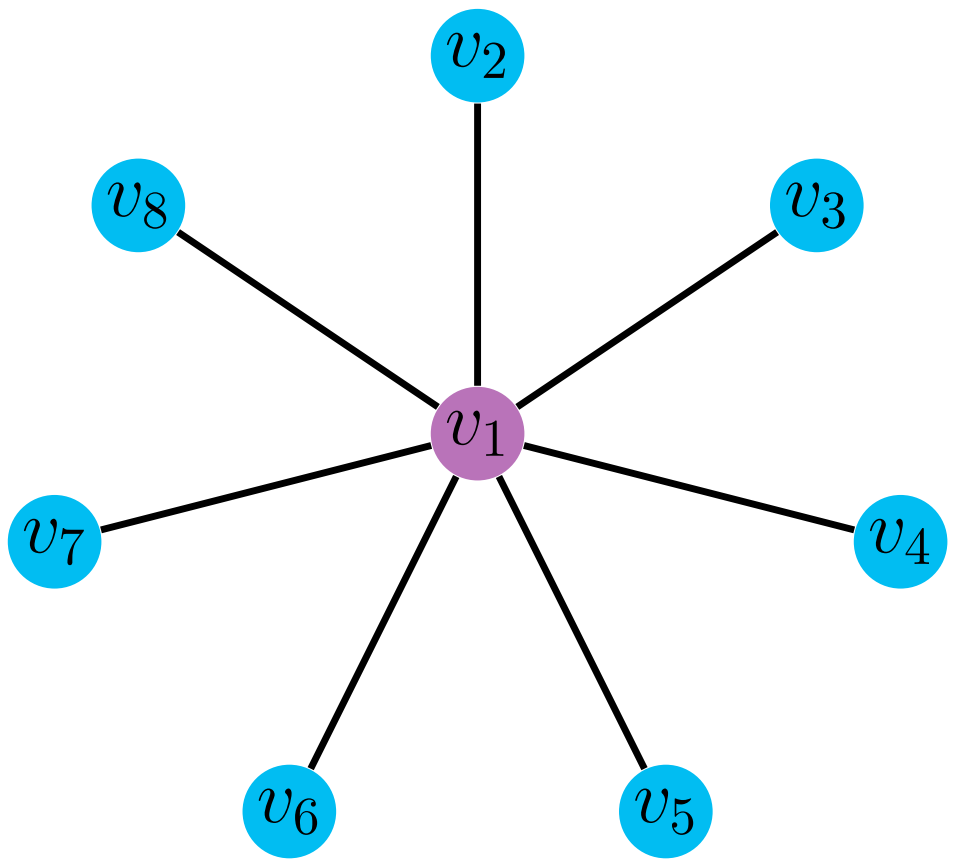

After adding edges connecting every pair of boundary vertices, the star $\stargraph_{n}$ becomes the complete graph $K_n$, so is indeed 3-connected in the sense of \cite{KenyonEnharmonic} for $n\geq4$. It is worth noting that considering star graphs is sufficient to establish that the action of $\Gal(\QQ^{\tr}/\QQ)$ is faithful; see \cite[Section 6]{KenyonEnharmonic}. Since the sum \eqref{eq:Lh} is taken over the vertices of $G$ adjacent to $v$, the discrete enharmonic Laplace equation for any graph with one interior vertex is the same as the discrete enharmonic Laplace equation for some star graph. Thus \Cref{prop:stargraphcounterexample} will show that any graph with a single interior vertex provides a particular counterexample to \Cref{conj:Conjecture1}. The star graph $\stargraph_8$ with seven leaves and eight vertices is shown in \Cref{fig:S8}.

\begin{proposition}\label{prop:stargraphcounterexample}
    For any star graph $\stargraph_{n}$ with $n\geq4$, the set of rational parameters with strictly positive edge energies and pairwise distinct boundary values for which the action of $\Gal(\QQ^{\tr}/\QQ)$ on the enharmonic functions is not transitive is Zariski dense in $\mathbb{A}_{\QQ}^{2n-2}$.
\end{proposition}
\begin{proof}
    Labeling the interior vertex as $v_1$ and the boundary vertices as $v_2,\dots, v_n$, the enharmonic functions of $\stargraph_{n}$ are the solutions of the single equation
\begin{equation}\label{eq:stargraph6polynomial}
        \sum_{j=2}^{n}\energy_{1j}\prod_{\substack{2\leq i\leq n\\i\neq j}} (x_1-b_i)=0
\end{equation}
where $b_j=b(v_j)$, $\energy_{1j}=\energy_{v_1v_j}$ and $x_1:=h(v_1)$.

We claim that the set of points whose coordinates give a rational root is Zariski dense in the parameter space. 
Consider the map
\begin{equation}\label{eq:phiStar}
    \begin{aligned}
        \varphi:\mathbb{A}_\QQ^{2n-2}&\dashrightarrow\mathbb{A}_\QQ^{2n-2}\\
        (b_2,b_3,\dots,b_n,y_1,\energy_{13},\dots,\energy_{1n})&\mapsto(b_2,b_3,\dots,b_n,\energy_{12},\energy_{13},\dots,\energy_{1n})
    \end{aligned}
\end{equation}
with 
\begin{equation}\label{eq:imageEnergy}
    \energy_{12}=-(y_1-b_2)\sum_{j=3}^n\frac{\energy_{1j}}{y_1-b_j}
\end{equation}
defined on the open subset where $y_1\neq b_2,b_3,\dots,b_n$.

Since there is no non-trivial polynomial relation among the pullbacks of the coordinates of the codomain, $\varphi$ is dominant.
A rational point in the domain is mapped to a point whose coordinates are the rational energies and boundary values giving the rational root $y_1$.

There is a Zariski-dense subset of the domain 
whose image has $\energy_{12}$ positive. Consider the subset of $\QQ^{2n-2}$ consisting of points such that
\begin{equation}\label{eq:inequalities}
b_2<y_1<b_3<\cdots<b_n
\qquad\text{and}\qquad
\energy_{1j}>0\quad (3\leq j\leq n).
\end{equation}
Equation \eqref{eq:imageEnergy} guarantees that the image of every point in this subset has $\energy_{12}>0$. It is Zariski dense in the domain because the inequalities \eqref{eq:inequalities} define a Euclidean-open subset of the real points of the domain in which the rational points are Euclidean, and therefore Zariski, dense.
Since $\varphi$ is dominant, we conclude that the set of points in $\QQ^{n-1}\times\QQ^{n-1}_{>0}$ such that \eqref{eq:stargraph6polynomial} has a rational root, and no root is equal to a boundary value, is Zariski dense in $\mathbb{A}_{\QQ}^{2n-2}$.

\end{proof}

When $G$ has two boundary vertices, the enharmonic functions are in bijection with the bipolar orientations of $G$ by \cite[Theorem 1]{KenyonEnharmonic}. Lutz pointed out \cite[Question 3.20]{lutz2019electrical} that if such a graph is also 3-connected, then \cite[Theorem 7.1]{deFraysseix} implies the connectedness of a graph obtained from the Dirichlet hyperplane arrangement $\mathcal{A}_{G,b}$ (\cite[Theorem 1.4]{lutz2019electrical}) by placing a vertex in every bounded chamber of $\mathcal{A}_{G,b}$ and adding an edge between any two vertices whose associated chambers share a facet. For this reason, such graphs make hopeful candidates for a guarantee of transitivity in terms of bipolar orientations and connectedness. However, an obstruction to the existence of a generic subset of the parameter space where the Galois action is transitive still exists in this case.

\begin{definition}
    The \demph{wheel graph} $W_{n}$ is the graph on $n$ vertices formed by connecting a universal vertex to each of the vertices of an $(n-1)$-cycle.
\end{definition}

In \cite[Proposition 7.2]{lutz2019electrical}, it was shown that if $G$ is a wheel graph on $n\equiv 3\pmod 4$ vertices with two boundary vertices on opposite sides of the rim of the wheel, then fixing all edge energies to one guarantees that the Galois action is not transitive. Such a graph can also be extended to a counterexample to \cite[Conjecture 1]{KenyonEnharmonic}.

\begin{remark}\label{rmk:WheelGraph}
    On the wheel graph $W_7$ (\Cref{fig:W7}), designate two opposite vertices on the outer $6$-cycle as boundary vertices and all other vertices as interior (\Cref{fig:W7}). The set of rational parameters with strictly positive edge energies and pairwise distinct boundary values for which the action of $\Gal(\QQ^{\tr}/\QQ)$
    is not transitive is Zariski dense in $\mathbb{A}_{\QQ}^{14}$.

    \begin{figure}
    \centering
    \input{W7}
    \caption{Wheel graph on seven vertices. Boundary vertices are marked in blue and interior vertices are marked in purple.}
    \label{fig:W7}
\end{figure}
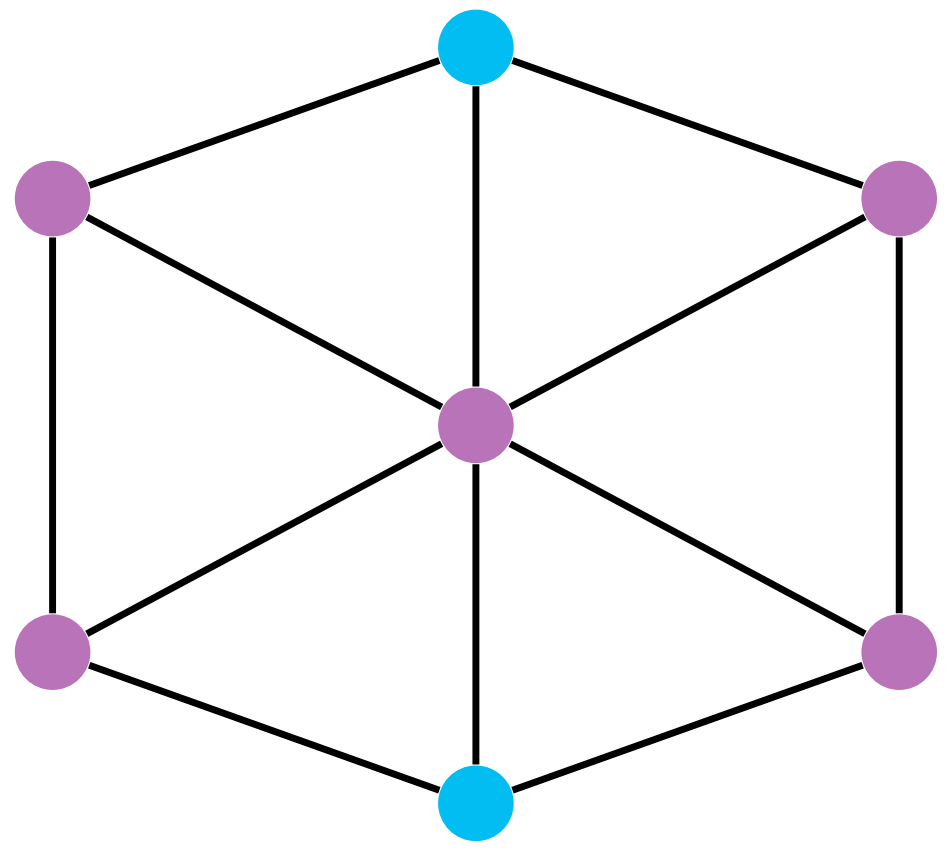
Similarly to \Cref{prop:stargraphcounterexample}, \Cref{rmk:WheelGraph} holds because there is a Zariski-dense subset of the parameter space such that an enharmonic function of $W_7$ has a coordinate with a rational root.

\end{remark}

\section{Generic Galois Groups}\label{sec:GenericGaloisGroups}
We refer to the Galois group of the discrete enharmonic Laplace equations of $G$ away from the thin set in \Cref{prop:3.3.5} as \demph{the Galois group of $G$}. It is isomorphic to the Galois group of the Galois closure of the extension of function fields given by $\pi$, also called the generic Galois group.
Structural features of the graph can force the generic Galois group to be a proper subgroup of the full symmetric group.

\subsection{Complete bipartite graphs}

\begin{proposition}
    Let $G=K_{m,n}$, $m\geq2,n\geq3$ be a complete bipartite graph such that the boundary vertices of $G$ all lie in one part and the interior vertices of $G$ lie in the complement. Then the Galois group of $G$ is not the full symmetric group.
\end{proposition}
\begin{proof}
 Label the part of size $m$ by $[m]=\{1,2,\dots,m\}$ and the part of size $n$ by $[n]$. If the part of size $n$ constitutes the boundary vertices, the discrete enharmonic Laplace equations for these complete bipartite graphs with boundary are
\begin{equation}\label{eq:K_{m,n}}
    \sum_{j=1}^n\frac{\energy_{ij}}{x_i-b_j}=0 \quad\text{for $i\in[m]$}.
\end{equation}
Since the equations \eqref{eq:K_{m,n}} are decoupled, the Galois group is a subgroup of a direct product of symmetric groups and so cannot be $S_k$ for $k\geq2$. 
\end{proof}
The Galois group may be a product of symmetric groups, as in the following example.
\begin{example}
    Let $G$ be the complete bipartite graph $K_{2,3}$ with the subset of size three constituting the boundary vertices. Then the Galois group of $G$ is the Klein four-group $V_4$.    
\end{example}
This example also shows that the Galois group of the enharmonic Laplace equations of a graph need not be at least alternating.
\subsection{The four-leaf star}\label{sec:TheFour-leafStar}

 We now illustrate a choice of thin set for a graph whose Galois group is the full symmetric group. For additional examples that illustrate finding thin sets in the context of Hilbert's irreducibility theorem, see \cite[Section 3.3]{serre}. 
 
\begin{example}
    For the star graph $\stargraph_5$, let $f$ be the monic cubic obtained from \eqref{eq:stargraph6polynomial} by dividing each coefficient by $\energy_{12}+\energy_{13}+\energy_{14}+\energy_{15}$ and let $\Delta_f$ be its discriminant. The generic Galois group of $\stargraph_5$ is the full symmetric group $S_3$. An explicit thin set $A\subset\mathbb{Q}^8$ can be taken as the union of the following loci:
    \begin{enumerate}[label=\textit{\roman*}.]
        \item the set where two boundary values coincide or an edge energy, the sum of the edge energies, or $\Delta_f$ vanishes;
    \item the set where $f_t$ has a rational root; and
    \item the set where $\Delta_f(t)$ is a square in $\mathbb{Q}$.
    \end{enumerate}
    The first locus is thin of type I and the latter two are thin of type II. They arise from covers of degree three and two, respectively. For every rational parameter choice $t\notin A$, the polynomial $f_t$ is irreducible and its discriminant is not square, hence its Galois group is $S_3.$ 
\end{example}

\section{Conclusion}\label{sec:Conclusion}
In this article, we considered the relationship between the arithmetic and combinatorics of the Dirichlet problem defined on a graph. When asking for combinatorial properties of a graph that correspond to orbits of the Galois action, one may hope that the choice of parameters is as irrelevant as possible so that the graph alone captures the arithmetic behavior of the solution set. We showed that if the property of interest is having a single orbit under the action of $\Gal(\Bar{\QQ}/\QQ)$, the appropriate interpretation of ``as irrelevant as possible" is as the complement of a thin set, rather than a Zariski-open set, in the parameter space. Examples showcased that parameters giving a rational solution to \eqref{eq:Lh} can obstruct the latter. These results apply to enharmonic functions on resistor networks \cite{KenyonEnharmonic,lutz2019electrical}, but may admit analogues in related applications.

Extending the boundary value function
to take values in $\mathbb{C}$ yields an interpretation of
the edge energies
as the reactive power in an AC circuit comprised
of inductors and capacitors  \cite[Section 4]{KenyonChromatic}. Imposing the condition 
\begin{equation}\label{eq:momentumConservation}
    \sum_{w\sim v}\energy_{vw}=0
\end{equation}
for every interior vertex $v$ promotes the set of energies, together with the solutions, to be the holomorphic quadratic differentials of $G$ \cite{KenyonChromatic}. When $G$ has exactly three boundary vertices, \eqref{eq:Lh} can be interpreted as the gauge-fixed CHY scattering equations \cite{CHY} with the edge energies considered Mandelstam invariants subject to momentum conservation \eqref{eq:momentumConservation}. In particular, the dependence of \eqref{eq:Lh} on $G$ corresponds to the case where some Mandelstam invariants are set to zero \cite{Early_Pfister_Sturmfels_2025,betti2025graphical}, an operation with interpretations in algebraic statistics \cite{gross2014maximum} and in physics \cite{Early_Pfister_Sturmfels_2025,Cao2024} that provides insight into the asymptotic behavior of the solutions \cite{agostini2023likelihood,Cachazo2020}. The problem of counting solutions to systems of the form \eqref{eq:Lh} has received substantial attention. It remains open to understand the arithmetic of the solution sets, and we hope that further work will uncover relationships to the robust combinatorics in these contexts.

\subsubsection*{Acknowledgments}

We thank Jose Israel Rodriguez for helpful conversations and feedback on previous drafts. We are grateful to Bernd Sturmfels and Rick Kenyon for encouraging this work and to Rick Kenyon for pointing out a connection to \cite{betti2025graphical}. 

{\footnotesize This material is based upon work supported by the National Science Foundation Graduate Research Fellowship Program under Grant No. 2137424. Any opinions, findings, and conclusions or recommendations expressed in this material are those of the authors and do not necessarily reflect the views of the National Science Foundation. Support was also provided by the Graduate School and the Office of the Vice Chancellor for Research at the University of Wisconsin-Madison with funding from the Wisconsin Alumni Research Foundation, as well as by National Science Foundation Grant No. 2510307}.

\bibliographystyle{siam-no-dash-title-color-links} 
\bibliography{refs} 

\end{document}

%% file: S8.tex
\begin{tikzpicture}[
  cyanpt/.style={circle, fill=cyan!75, minimum size=5mm, inner sep=1pt},
  purppt/.style={circle, fill=violet!55, minimum size=5mm, inner sep=1pt},
  edge/.style={line width=1.3pt},
  lab/.style={font=\large}
]

\node[purppt, lab] (v1) at (0,0) {$v_1$};

\foreach \i [evaluate=\i as \j using int(\i+1)] in {1,...,7} {
  \node[cyanpt, lab] (v\j) at
    ({2.87*cos(90 - 360*(\i-1)/7)},
     {-0.13 + 2.63*sin(90 - 360*(\i-1)/7)})
    {$v_{\j}$};
}

\foreach \j in {2,...,8} {
  \draw[edge] (v1) -- (v\j);
}

\end{tikzpicture}

%% file: W7.tex
\begin{tikzpicture}[
  cyanpt/.style={circle, fill=cyan!75, minimum size=5mm, inner sep=0pt},
  purppt/.style={circle, fill=violet!55, minimum size=5mm, inner sep=0pt},
  edge/.style={line width=1.3pt},
  lab/.style={font=\large}
]

\node[cyanpt, label={[lab]below:}] (w1) at (0,-2.5) {};
\node[cyanpt, label={[lab]above:}] (w2) at (0,2.5) {};

\node[purppt, label={[lab]above left:}] (v7) at (-2.8,1.5) {};
\node[purppt, label={[lab]below left:}] (v6) at (-2.8,-1.5) {};
\node[purppt, label={[lab]above right:}] (v4) at (2.8,1.5) {};
\node[purppt, label={[lab]right:}] (v5) at (2.8,-1.5) {};
\node[purppt, label={[lab]above left:}] (v3) at (0,0) {};

\draw[edge] (w2)--(v7)--(v6)--(w1)--(v5)--(v4)--(w2);

\draw[edge] (v3)--(w2);
\draw[edge] (v3)--(w1);
\draw[edge] (v3)--(v7);
\draw[edge] (v3)--(v6);
\draw[edge] (v3)--(v4);
\draw[edge] (v3)--(v5);

\end{tikzpicture}